\documentclass{gtart_h}  

\def\ifplaintex{\expandafter\ifx\csname documentclass\endcsname\relax}

\ifplaintex 
\else
\fi

\expandafter\ifx\csname beginpicture\endcsname\relax
\expandafter\ifx\csname documentclass\endcsname\relax
\input pictex \else
\input prepictex \input pictex \input postpictex \fi\fi

\def\gt{{\mathsurround=0pt\it $\cal G\mskip-2mu$eometry \&\ 
$\cal T\!\!$opology}}        

\def\gtp{{\mathsurround=0pt\it $\cal G\mskip-2mu$eometry \&\ 
$\cal T\!\!$opology $\cal P\!$ublications}}  

\def\lognumber#1{\def\thelognumber{#1}}
\def\volumenumber#1{\def\thevolumenumber{#1}}
\def\papernumber#1{\def\thepapernumber{#1}}
\def\volumeyear#1{\def\thevolumeyear{#1}}

\def\pagenumbers#1#2{\def\startpage{#1}\def\finishpage{#2}}
\def\published#1{\def\publishdate{#1}}
\def\proposed#1{\def\theproposer{#1}}
\def\seconded#1{\def\theseconders{#1}}
\def\received#1{\def\receiveddate{#1}}
\def\revised#1{\def\reviseddate{#1}}
\def\accepted#1{\def\accepteddate{#1}}
\def\asciititle#1{\def\theasciititle{#1}}
\def\covertitle#1{\def\thecovertitle{#1}}

\long\def\asciiabstract#1{\long\def\theasciiabstract{#1}}
\def\asciikeywords#1{\def\theasciikeywords{#1}}

\let\\\par\let\thelognumber\relax
\let\thevolumenumber\relax\let\thepapernumber\relax
\let\thevolumeyear\relax\let\thesamplenumber\relax\let\startpage\relax
\let\finishpage\relax\let\publishdate\relax\let\receiveddate\relax
\let\reviseddate\relax\let\accepteddate\relax\let\theasciititle\relax
\let\thecovertitle\relax\let\theasciiauthors\relax
\let\theasciiabstract\relax\let\theasciikeywords\relax
\let\theasciiemail\relax\let\theshortauthors\relax\let\theshorttitle\relax

\long\def\maketitlep{   

\count0=\startpage

\gt\hfill      
\beginpicture
\setcoordinatesystem units <0.33truein, 0.33truein> point at 2.2 0.9
\setplotsymbol ({$\cal G$})
\plotsymbolspacing=9truept
\circulararc 315 degrees from 0 1 center at 0 0
\setplotsymbol ({$\cal T$})
\circulararc 315 degrees from 1 -1 center at 1 0
\endpicture
\break
{\small\ifx\thesamplenumber\relax 
Volume \else Sample
\fi\thevolumenumber\ (\thevolumeyear)
\startpage--\finishpage\nl
Published: \publishdate}
\vglue 0.5truein plus 0.4fil minus 0.1truein

{\parskip=0pt\leftskip 0pt plus 1fil\def\\{\par\smallskip}{\ifplaintex\large
\else\Large\fi\bf\thetitle}\par\medskip}   

\vglue 0pt plus 0.1fil 

{\parskip=0pt\leftskip 0pt plus 1fil\def\\{\par}{\sc\theauthors}
\par\medskip}

\vglue 0pt plus 0.1fil 

{\small\parskip=0pt\let\newline\\
{\leftskip 0pt plus 1fil\def\\{\par}{\sl\theaddress}\par}
\expandafter\ifx\theemail\relax    
\relax\else\vglue 5pt plus 0.02fil minus 2pt\def\\{\stdspace{\rm 
and}\stdspace} 
\cl{Email:\stdspace\tt\theemail}\fi
\ifx\theurl\relax                  
\relax\else\vglue 5pt plus 0.02fil minus 2pt\def\\{\stdspace{\rm 
and}\stdspace}
\cl{URL:\stdspace\tt\theurl}\fi\par}

\vglue 7pt plus 0.3fil minus 3pt

{\bf Abstract}
\vglue 5pt plus 0.1fil minus 2pt

\theabstract

\vglue 7pt plus 0.3fil minus 3pt

{\bf AMS Classification numbers}\quad Primary:\quad \theprimaryclass

Secondary:\quad \thesecondaryclass

\vglue 5pt plus 0.3fil minus 2pt

{\bf Keywords:}\quad \thekeywords

\vglue 10pt plus 0.5fil minus 5pt

{\small  Proposed: \theproposer\hfill Received: \receiveddate\nl
Seconded: \theseconders\hfill 
\ifx\reviseddate\relax                         
Accepted: \accepteddate                        
\else
Revised: \reviseddate                          
\fi}
\eject
}       

\let\maketitlepage\maketitlep
\let\maketitle\maketitlepage

\font\phead=cmsl9 scaled 950
\font=cmsl9 scaled 1050
\font\pnum=cmbx10 scaled 913
\font\lnum=cmbx10 
\font\pfoot=cmsl9 scaled 950
\font=cmsl9 scaled 1050
\ifplaintex
\headline{\vbox to 0pt{\vskip -4.5mm\line{\small\phead\ifnum
\count0=\startpage ISSN 1364-0380 (on line)
1465-3060 (printed) \hfill {\pnum\folio}\else\ifodd\count0\def\\{ }%
\ifx\theshorttitle\relax\thetitle\else\theshorttitle\fi\hfill{\pnum\folio}
\else\def\\{ and }{\pnum\folio}\hfill\ifx\theshortauthors\relax\theauthors
\else\theshortauthors\fi\fi\fi}\vss}}
\footline{\vbox to 0pt{\vglue 0mm\line{\small\pfoot\ifnum\count0=\startpage
\copyright\ \gtp\hfill\else
\gt, Volume \thevolumenumber\ (\thevolumeyear)\hfill\fi}\vss
}}
\else
\makeatletter
\def\@oddhead{{\small\lhead\ifnum\count0=\startpage ISSN 1364-0380 (on line)
1465-3060 (printed) \hfill {\lnum\number\count0}\else\ifodd\count0
\def\\{ }\ifx\theshorttitle\relax \thetitle \else\theshorttitle\fi\hfill
{\lnum\number\count0}\else\def\\{ and }{\lnum\number\count0}
\hfill\ifx\theshortauthors\relax 
\theauthors\else\theshortauthors\fi\fi\fi}}\def\@evenhead{\@oddhead}
\def\@oddfoot{\small\lfoot\ifnum\count0=\startpage\copyright\ \gtp\hfill\else
\gt, Volume \thevolumenumber\ (\thevolumeyear)\hfill\fi}
\def\@evenfoot{\@oddfoot}
\makeatother
\fi

\newwrite\gtoutfile
\long\gdef\makeheadfile{  
{\def\\{, }\def\s{ }
\immediate\openout\gtoutfile head.xxx
\immediate\write\gtoutfile{Proxy-for: \ifx\theasciiauthors\relax
\theauthors\else\theasciiauthors\fi\s<\ifx\theasciiemail\relax\theemail\else\theasciiemail\fi>}
\immediate\write\gtoutfile{\noexpand\\}
\immediate\write\gtoutfile{Authors: \ifx\theasciiauthors\relax
\theauthors\else\theasciiauthors\fi}
{\def\\{ }\immediate\write\gtoutfile{Title: \ifx\theasciititle\relax
\thetitle\else\theasciititle\fi}}
\immediate\write\gtoutfile{Subj-class: GT or SG or MG etc}
\immediate\write\gtoutfile{MSC-class: \theprimaryclass\ifx\thesecondaryclass\relax\else, \thesecondaryclass\fi}
\immediate\write\gtoutfile{Journal-ref: Geom. Topol. \thevolumenumber
(\thevolumeyear) \startpage-\finishpage}
\immediate\write\gtoutfile{Comments: Published by Geometry and Topology at}
\immediate\write\gtoutfile{\s\s http://www.maths.warwick.ac.uk/gt/GTVol\thevolumenumber/paper\thepapernumber.abs.html}
\immediate\write\gtoutfile{\noexpand\\}
\immediate\write\gtoutfile{}
\ifx\theasciiabstract\relax
\immediate\write\gtoutfile{\theabstract}\else
\immediate\write\gtoutfile{\theasciiabstract}\fi
\immediate\write\gtoutfile{}
\immediate\write\gtoutfile{\noexpand\\}
\immediate\write\gtoutfile{}
\immediate\closeout\gtoutfile}}  

\def\maketitlepage{\maketitlep\makeheadfile}
\let\maketitle\maketitlepage

\lognumber{432}

\volumenumber{8}\papernumber{28}\volumeyear{2004}
\pagenumbers{1032}{1042}
\received{15 March 2004}
\published{7 August 2004}
\revised{3 August 2004}
\accepted{11 July 2004}
\proposed{Benson Farb}
\seconded{Walter Neumann, Steven Ferry}

\usepackage{amsmath, amsfonts}

\renewcommand\theenumi{\arabic{enumi}}

\DeclareMathOperator{\I}{\bf I}
\DeclareMathOperator{\II}{\bf II}
\DeclareMathOperator{\III}{\bf III}

\DeclareMathOperator{\IV}{\bf IV}

\DeclareMathOperator{\V}{\bf V}

\DeclareMathOperator{\C}{\bf C}

\DeclareMathOperator{\im}{Im}
\DeclareMathOperator{\Ker}{Ker}

\DeclareMathOperator{\Stab}{Stab}
\DeclareMathOperator{\GHS}{GHS}

\def\b{\operatorname{b}}

\newcommand\BS{\mathbb{S}}

\newcommand\cn{\mathcal{N}_\bq}
\newcommand\Sym[1]{\mathbf{S}_{#1}}

\newtheorem{Theorem}[subsection]{Theorem}
\newtheorem{Prop}[subsection]{Proposition}

\newtheorem{Lemma}[subsection]{Lemma}

\newtheorem*{I(n)}{$\I(n)$}
\newtheorem*{II(n)}{$\II(n)$}
\newtheorem*{III(2k+1)}{$\III(2k+1)$}
\newtheorem*{III'(2k+1)}{$\III'(2k+1)$}
\newtheorem*{III''(2k+1)}{$\III''(2k+1)$}

\newtheorem*{IV(n)}{$\IV(n)$}
\newtheorem*{IV(2k+1)}{$\IV(2k+1)$}
\newtheorem*{V(n)}{$\V(n)$}
\newtheorem*{C(n)}{$\C(n)$}

\newcommand\bq{{\mathbf q}}
\newcommand\bone{{\mathbf 1}}
\newcommand\Ht{\operatorname{\mathfrak h}^\bq}

\def\H{\operatorname{\mathfrak h}}

\def\h{\operatorname{\ltwo{\mathcal H}}}
\def\hq{\operatorname{\ltwoq\mathcal H}}
\def\bt{\operatorname{b}_\bq}
\newcommand{\ltwo}{{\mathit L^2}}
\newcommand{\ltwoq}{{\mathit L^2_\bq}}
\newcommand{\GSC}{Generalized Singer Conjecture}
\newcommand{\SC}{Singer Conjecture}

\theoremstyle{definition}

\newtheorem*{Hopf Conjecture}{Hopf Conjecture}
\newtheorem*{The Flag Complex Conjecture}{The Flag
Complex Conjecture}
\newtheorem*{SingerConjecture}{The \SC}
\newtheorem*{GS}{The \GSC}
\newtheorem*{Atiyah Conjecture}{Atiyah Conjecture}
\newcommand{\bdd}[1]{{b\partial \Delta^{#1}}}
\newcommand{\bddv}[1]{{b\partial \Delta_v^{#1}}}

\def\l{\operatorname{\ell}}

\title{Weighted $\ltwo$--cohomology of Coxeter groups\\based 
on barycentric subdivisons }
\asciititle{Weighted L^2-cohomology of Coxeter groups
based on barycentric subdivisons}
\covertitle{Weighted $L^2$--cohomology of Coxeter groups
based on barycentric subdivisons}

\author{Boris L Okun}
\address{Department of Mathematical Sciences\\
        University of Wisconsin--Milwaukee\\
        Milwaukee, WI 53201, USA}
\email{okun@uwm.edu}



\asciikeywords{Coxeter group, aspherical manifold, barycentric
subdivision, weighted L^2-cohomology, Tomei manifold, Singer
conjecture}

\begin{abstract}
Associated to any finite flag complex $L$ there is a right-angled
Coxeter group $W_L$ and a contractible cubical complex $\Sigma_L$
(the Davis complex) on which $W_L$ acts properly and cocompactly,
and such that the link of each vertex is $L$. It follows that if
$L$ is a generalized homology sphere, then $\Sigma_L$ is a
contractible homology manifold. We prove a generalized version of
the Singer Conjecture (on the vanishing of the reduced weighted
$\ltwoq$--cohomology above the middle dimension) for the
right-angled Coxeter groups based on barycentric subdivisions  in
even dimensions. We also prove this conjecture for the groups
based on the barycentric subdivision of the boundary complex of a
simplex.
\end{abstract}

\asciiabstract{Associated to any finite flag complex L there is a
right-angled Coxeter group W_L and a contractible cubical complex
Sigma_L (the Davis complex) on which W_L acts properly and
cocompactly, and such that the link of each vertex is L. It follows
that if L is a generalized homology sphere, then Sigma_L is a
contractible homology manifold.  We prove a generalized version of the
Singer Conjecture (on the vanishing of the reduced weighted
L^2_q-cohomology above the middle dimension) for the right-angled
Coxeter groups based on barycentric subdivisions in even
dimensions. We also prove this conjecture for the groups based on the
barycentric subdivision of the boundary complex of a simplex.}
\begin{document}

\maketitlepage

\section{Introduction}

A construction of Davis (\cite{D1}, \cite{D2}, \cite{D3}),
associates to any finite flag complex $L$,   a ``right-angled''
Coxeter group $W_L$ and a contractible cubical cell complex
$\Sigma_L$ on which $W_L$ acts properly and cocompactly. $W_L$ has
the following presentation: the generators are the vertices of
$L$, each generator has order $2$, and two generators commute if
the span an edge in $L$. The  most important feature of this
construction is that the link of each vertex of $\Sigma_L$ is
isomorphic to $L$. A simplicial complex $L$ is a \emph{generalized
homology $m$--sphere} (for short, a $\GHS^m$) if it is a homology
$m$--manifold having the same homology as a standard sphere $\BS^m$
(the homology is with real coefficients.) It follows, that if $L$
is a $\GHS^{n-1}$, then $\Sigma_L$ is a homology $n$--manifold.

If $L$ is a simplicial complex,  $bL$ will denote the barycentric
subdivision of $L$. $bL$ is a flag simplicial complex.  Let
$\partial\Delta^n$ denote the boundary complex of the standard
$n$--dimensional simplex.

We study a certain  weighted $\ltwo$--cohomology theory $\hq^*$,
described in \cite{Dy}, \cite{DDJO}. Suppose, for each vertex of
$v\in L$ we are given a positive real number $q_v$, and let $\bq$
denote the vector with components $q_v$. Given a minimal word
$w=v_1 \dots v_n\in W_L$, let $\bq^w$=$q_{v_1} \dots q_{v_n}$. For
each $W_L$--orbit of cubes pick a representative $\sigma_0$ and let
$w(\sigma)=w$ if $\sigma=w\sigma_0$. (The ambiguity in the choices
will not matter in our discussion.) Let $\ltwoq
C^i(\Sigma_L)=\{\Sigma c_\sigma \sigma \mid \Sigma c^2_\sigma
\bq^{w(\sigma)} < \infty\}$ be the Hilbert space of infinite
$i$--cochains, which are square-summable with respect to the weight
$\bq^w$. The usual coboundary operator $d$ is then a bounded
operator, and we define the reduced weighted $\ltwoq$--cohomology
to be $\hq^i(\Sigma_L)=\Ker(d^i)/\overline{\im(d^{i-1})}$.
Similarly, one can define  the reduced weighted $\ltwoq$--homology,
except, instead of the usual boundary operator one uses the
adjoint of $d$. It follows from the Hodge decomposition that the
resulting homology and cohomology spaces are naturally isomorphic.
These spaces are Hilbert modules over the Hecke--von Neumann
algebra $\cn$ (an appropriately completed Hecke algebra of $W_L$.)
This allows us to introduce the weighted $\ltwoq$ Betti numbers
--- the dimension of $\hq^i$ over $\cn$. If
$\bq=\bone=(1,\dots,1)$, we obtain the usual reduced
$\ltwo$--cohomology, and we omit the index $\bq$. We write $\bq \le
\bone$, if each component of $\bq$ is $\le 1$.

The following conjecture, attributed to Singer, goes back to
1970's.
\begin{SingerConjecture}
If $M^n$ is a closed aspherical manifold, then
$$\h^i(\widetilde{M}^n)=0 \text{ for all $i\not= n/2$.}$$
\end{SingerConjecture}

As explained in \cite[Section 14]{DDJO}, the appropriate
generalization of the Singer Conjecture to the weighted case is
the following conjecture:

\begin{GS}\label{I}
Suppose $L$ is a flag $\GHS^{n-1}$. Then $\hq^i(\Sigma_L)=0$ for
$i> n/2$ and $\bq \le \bone$.
\end{GS}
(Poincar\'e duality shows that for $\bq=\bone$ this conjecture
implies the Singer Conjecture for $\Sigma$.)

This conjecture holds true for $n\le 4$ by \cite{DDJO}. One of the
main results of this paper is a proof of this conjecture for
barycentric subdivisions in even dimensions. The proof uses a
reduction to a very special case $L=\bdd{2k-1}$.

We prove this case as Theorem \ref{M}. It turns out (Theorem
\ref{III}), that this result implies  the vanishing of the
$\ltwoq$--cohomology in a certain range for arbitrary right-angled
Coxeter groups based on barycentric subdivisions. (For
$\bq=\bone$, this implication is proved in \cite{DO2}.) In
particular, it follows that the \GSC\ is true for all barycentric
subdivisions in even dimensions (Theorem \ref{IV}), and  for
$\bdd{n}$ in all dimensions (Theorem \ref{V}).

This paper relies very heavily on \cite{DDJO}. In the inductive
proofs we mostly omit the first steps, they are easy exercises
using \cite{DDJO}.

The author wishes
to thank Mike\,Davis, Jan\,Dymara and Tadeusz\,Januszkiewicz for
several illuminating discussions.

\section{Vanishing conjectures}

We will follow the notation from \cite{DDJO}. Given a flag complex
$L$ and a full subcomplex $A$, set:
\begin{align*}
\Ht_i(L) &= \h_i(\Sigma_L)\\
\Ht_i(A) &= \h_i(W_L\Sigma_A)\\
\Ht_i(L,A) &= \h_i(\Sigma_L, W_L\Sigma_A)\\
\bt^i(L) &= \dim_{\cn}(\Ht_i(L))\\
\bt^i(L,A) &= \dim_{\cn}(\Ht_i(L,A))
\end{align*}

The dimension of $\Sigma_L$ is one greater than the dimension of
$L$. Hence, $\bt^i(L)=0$ for $i>\dim L +1$.

We will use the following three properties of $\ltwoq$--homology.

\begin{Prop}[{See \cite[Section 15]{DDJO}}]
\noindent
\begin{description}
\item[The Mayer--Vietoris sequence]\label{MVS}
If $L=L_1\cup L_2$ and $A=L_1\cap L_2$, where $L_1$ and  $L_2$
(and therefore, $A$) are full subcomplexes of $L$, then
$$\to\Ht_i(A)\to\Ht_i(L_1)\oplus\Ht_i(L_2)\to\Ht_i(L)\to$$
is weakly exact.

\item[The K\"unneth Formula]\label{KF}
The Betti numbers of the join of two complexes are given by:
$$\bt^k(L_1\ast L_2) = \sum_{i+j=k}\bt^i(L_1)\bt^j(L_2) .$$

\item[Poincar\'e Duality]\label{PD}
If $L$ is a flag $\GHS^{n-1}$, then
$\bt^i(L)=\b_{\bq^{-1}}^{n-i}(L)$.
\end{description}
\end{Prop}

If $\sigma$ is a simplex in $L$, let $L_\sigma$ denote the link of
$\sigma$ in $L$. To simplify notation we will write $bL_v$ instead
of $(bL)_v$ to denote the link of the vertex $v$ in $bL$. Let
$\mathcal{C}$ be a class of $\GHS$'s closed under the operation of
taking link of vertices, i.e.\ if $S\in\mathcal{C}$ and $v$ is a
vertex of $S$ then $S_v\in\mathcal{C}$. Following Section 15 of
\cite{DDJO} we consider several variations of the \GSC\ for the
class $\mathcal{C}$.

\begin{I(n)}
If $S\in \mathcal{C}$ and $\dim S=n-1$, then  $\bt^i(S)=0$ for
$i>n/2$ and $\bq \le \bone$.
\end{I(n)}

\begin{III'(2k+1)}
Let $S\in \mathcal{C}$ and $\dim S=2k$. Let $v$ be a vertex of
$S$. Then the map $i_\ast\co\Ht_k(S_v)\to \Ht_k(S)$, induced by
the inclusion, is the zero homomorphism for $\bq \ge \bone$.
\end{III'(2k+1)}

\begin{V(n)}
Let $S\in \mathcal{C}$ and $\dim S=n-1$. Let $A$ be a full
subcomplex of $S$.
\begin{itemize}
\item If $n=2k$ is even, then $\bt^i(S,A)=0$ for all $i>k$
and $\bq \le \bone$.
\item If $n=2k+1$ is odd, then $\bt^i(A)=0$ for all $i>k$
and $\bq \le \bone$.
\end{itemize}
\end{V(n)}

The argument in Section 16 of \cite{DDJO} goes through without
changes if we consider only $\GHS$'s from a class $\mathcal{C}$ to
give the following:

\begin{Theorem}[{Compare \cite[Section 16]{DDJO}}] \label{t31}
If we only consider
$\GHS$'s from a class $\mathcal{C}$, then the following
implications hold.
\begin{enumerate}
\item
$\I(2k+1)\implies \III'(2k+1 )$.
\item
$\V(n)\implies \I(n)$.
\item
$\V(2k-1)\implies \V(2k)$.
\item
$\left[\V(2k)\text{ and }\III'(2k+1)\right]\implies \V(2k+1) .$
\end{enumerate}
\end{Theorem}

Let $\mathcal{JD}$ denote the class of finite joins of the
barycentric subdivisions of the boundary complexes of standard
simplices:
$$\mathcal{JD}=\{\bdd{n_1}*\dots*\bdd{n_i}\}.$$

\begin{Lemma}\label{l:2.3}
The class $\mathcal{JD}$ is closed under the operation of taking
link of vertices.
\end{Lemma}

\begin{proof}
Let $S=\bdd{n_1}*\dots*\bdd{n_j}$ and $v \in S$.  We can assume
that $v\in \bdd{n_1}$. Then
$S_v=\bddv{n_1}*\bdd{n_2}*\dots*\bdd{n_j}=
\bdd{\dim(v)}*\bdd{n_1-\dim(v)-1}*\bdd{n_2}*\dots*\bdd{n_j}$, and
therefore $S\in \mathcal{JD}$.
\end{proof}

Next, consider the following statement:

\begin{III''(2k+1)}
Let $v$ be a vertex of $\bdd{2k+1}$. Then the map
$i_\ast\co\Ht_k(\bddv{2k+1})\to \Ht_k(\bdd{2k+1})$, induced by the
inclusion, is the zero homomorphism for $\bq \ge \bone$.
\end{III''(2k+1)}

\begin{Lemma}\label{l:2.4}
$\III''(2k+1)\implies \III'(2k+1)$ for the class$\mathcal{JD}$.
\end{Lemma}

\begin{proof}
By induction, we can assume that the lemma holds for all odd
numbers $<2k+1$. Then it follows from the Theorem~\ref{t31} that
$\V(m)$ and therefore $\I(m)$ hold for all $m<2k+1$.

Let $S=\bdd{n_1}*\dots*\bdd{n_j}$ with $n_1+\dots+n_j=2k+1$ and $v
\in S$. We assume that $v\in \bdd{n_1}$. Then
$S_v=\bddv{n_1}*\bdd{n_2}*\dots*\bdd{n_j}$  and, by the
K{\"u}nneth formula, the map in question  decomposes as the direct
sum of maps of the form
$$ (\Ht_{k_1}(\bddv{n_1})\to \Ht_{k_1}(\bdd{n_1}))
\otimes \bigotimes_{i=2}^j(\Ht_{k_i}(\bdd{n_i})\to
\Ht_{k_i}(\bdd{n_i}))$$
where $k_1+\dots+k_j=k$. Since $n_1+\dots+n_j=2k+1$ it follows
that $k_i<n_i/2$ for some index $i$. If $n_i<2k+1$, then the range
of the corresponding map in the above tensor product is $0$ by
$\I(n_i)$ and Poincar\'{e} duality, and therefore the tensor
product map is $0$. If $n_i=2k+1$ then, in fact, $i=1$ (the join
is a trivial join) and the result follows from $\III''(2k+1)$.
\end{proof}

Thus, it follows from Theorem \ref{t31}, Lemmas \ref{l:2.3} and
\ref{l:2.4}, and induction on dimension, that in order to prove
the \GSC\ for the class $\mathcal{JD}$ all we need is to prove
$\III''(2k+1)$.

\section{Removal of an odd-dimensional vertex}

Let $L$ be a simplicial complex and $bL$ be its barycentric
subdivision. The vertices of $bL$ are naturally graded by
"dimension": each vertex $v$ of $bL$ is the barycenter of a unique
cell (which we still denote $v$) of the complex $L$, and we call
the dimension of this cell the \emph{dimension} of the vertex $v$.
Let $E_L$ denote the subcomplex of $bL$ spanned by the even
dimensional vertices. Let $\mathcal{A}_L$ denote the set of full
subcomplexes $A$ of $L$ containing $E_L$, which have the following
property: if $A$ contains a vertex of odd dimension $2j+1$, then
$A$ contains all vertices of $bL$ of dimensions $\le 2j$. In other
words, any such $A$ can be obtained inductively from $bL$ by
repeated removal of an odd-dimensional vertex of the highest
dimension.

If $L=\partial\Delta^{n}$ we will use the notation   $E_n=E_L$ and
$\mathcal{A}_n =\mathcal{A}_L$.

\begin{Lemma}\label{An}
Assume $\III''(2m+1)$ holds  for $2m+1<n$.  Then for any
$(n-1)$--dimensional simplicial complex $L$ and any complex
$A\in\mathcal{A}_L$ we have:
 $$\bt^i(A)=\bt^i(bL)=0 \text{ for $i> (n+1)/2 $ and $\bq \le \bone$.}$$
\end{Lemma}

\begin{proof}
By induction, we can assume that the lemma holds for all $m<n$.
First, we claim that removal of odd-dimensional vertices does not
change the homology above $ (n+1)/2 $. Let $A\in\mathcal{A}_{L}$
and let $B=A-v$ where $v$ is a vertex of the highest odd dimension
of $A$. We let $\dim(v)=2d-1$, $1\le d \le k$.  We want to prove
that $\bt^i(A)=\bt^i(B)$ for $i> (n+1)/2  $. Consider the
Mayer--Vietoris sequence of the union $A=B\cup_{A_v} CA_v$:
$$\to\Ht_i(A_v)\to\Ht_i(B)\oplus\Ht_i(CA_v)\to\Ht_i(A)\to \Ht_{i-1}(A_v)$$
Suppose $i>  (n+1)/2 $.  Since $A_v=B\cap bL_v=B\cap
(\bdd{2d-1}*b(L_{_v}))$, and since $B\in \mathcal{A}_{L}$, it
follows, by construction, that $A_v$ splits as the join:
$$A_v=\bdd{2d-1}*A_1,$$ with
$A_1\in\mathcal{A}_{(L_v)}$. By  inductive assumption the lemma
holds for $L_v$, i.e.\ $\bt^i(A_1)=\bt^i(b(L_{_v}))= 0$ for $i>
(n+1)/2 -d$.

Since $\III''(2d-1)$ holds by hypothesis, by Lemma~\ref{l:2.4} and
Theorem \ref{t31}, $\I(2d-1)$ holds for the class $\mathcal{JD}$,
and, thus, $\bt^i(\bdd{2d-1})=0$ for $i\ge d$.

Then, by the K\"unneth formula, $\bt^{i-1}(A_v)=0$ for $i-1\ge
(n+1)/2 $, i.e.\ for $i> (n+1)/2 $. By \cite[Proposition
15.2(d)]{DDJO}, $\bt^i(CA_v)=\frac{1}{q_v+1}\bt^{i}(A_v)$.
Therefore in the above sequence the terms corresponding to $A_v$
and $CA_v$ are $0$, and the claim follows. Then it follows by
induction, that  $\bt^i(A)=\bt^i(bL)$ for all
$A\in\mathcal{A}_{L}$ and $i> (n+1)/2 $.

To prove the vanishing we note that, in particular,
$\bt^i(E_L)=\bt^i(bL)$ for $i> (n+1)/2 $. Since $E_L$ is spanned
by the even-dimensional vertices of $bL$ and since a simplex in
$bL$ has vertices of pairwise different dimensions, we have
$\dim(E_L)= \left[(n+1)/2\right]-1$. Therefore, $\bt^i(E_L)=0$ for
$i> (n+1)/2 $ and we have proved the lemma.
\end{proof}

In the special case $L=\Delta^{2k+1}$ this lemma admits the
following strengthening:

\begin{Lemma}\label{A2k+1}
Let $n=2k+1$. Assume $\III''(2m+1)$ holds  for $2m+1<n$. Then for
any complex $A\in\mathcal{A}_n$, $A\subset \bdd{n}$,   we have:
$$\bt^i(A)=\bt^i(\bdd{n}) \text{ for $i>k$ and $\bq \le \bone$.}$$
\end{Lemma}

\begin{proof}
We proceed as in the previous proof. As before, we have $B=A-v$,
$\dim(v)=2d-1$ and $A_v=\bdd{2d-1}*A_1$, where now $A_1\subset
\bdd{2k+1-2d}$. Therefore, the inductive assumption and the
hypothesis on $\III''(2d-1)$ imply that $b_i(A_1)=0$ for $i>k+d$.
The lemma follows as before.
\end{proof}

As explained in~\cite{DO2}, when $\bq = \bone$, the removal of the
odd-dimensional vertex does not change homology in \emph{all}
dimensions. We record this result below.

\begin{Lemma}\label{An1}
Assume $\III''(2m+1)$ holds  for $2m+1<n$ and $\bq = \bone$. Then
for any $(n-1)$--dimensional simplicial complex $L$ and for any
complex $A\in \mathcal{A}_L$, obtained by the repeated removal of
highest odd-dimensional vertices, we have:
 $$\b^*(A)=\b^*(bL).$$
\end{Lemma}

\begin{proof}
Again we repeat the proof of Lemma \ref{An}. As before, we have
the splitting $A_v=\bdd{2d-1}*A_1$. The point now is that for $\bq
= \bone$, $\I(2d-1)$ and Poincar\'{e} duality imply
$\b^*(\bdd{2d-1})=0$ and therefore $\b^*(A_v)=0$ by the K\"unneth
formula.
\end{proof}

\section{Intersection form}
\begin{Lemma}\label{ISO}
Let $L$ be a $\GHS^{2k}$ and let $v$ be a vertex of $L$. Then the
image of the restriction map on  $\ltwo$--cohomology $i^*\co
\h^k(\Sigma_L) \to \h^k(\Sigma_{L_v})$ is an isotropic subspace of
the intersection form of $\Sigma_{L_v}$.
\end{Lemma}

\begin{proof}
Note that the cup product of two $\ltwo$--cocycles is an
$L^1$--cocycle.  The intersection form is the result of evaluation
of the cup product of two middle-dimensional cocycles on the
fundamental class, which is $L^\infty$. Since $\Sigma_{L_v}$
bounds a half-space in $\Sigma_L$, $i_*([\Sigma_{L_v}])=0$ in
$L^\infty$--homology of $\Sigma_L$. Thus, if $a,b\in \h^n(\Sigma_L)
$, then $\langle i^*(a) \cup i^*(b), [\Sigma_{L_v}]\rangle =
\langle a\cup b, i_*([\Sigma_{L_v}])\rangle =0$.
\end{proof}

\begin{Lemma}\label{NIS}
Let $G$ be a group and let $A$ be a bounded $G$--invariant (with
respect to the diagonal action) non-degenerate bilinear form on a
Hilbert submodule $M \subset \l^2(G)$. Then $A$ has no nontrivial
$G$--invariant isotropic subspaces.
\end{Lemma}

\begin{proof}
Let us consider the case $M = \l^2(G)$ first.  $G$--invariance and
continuity of the form $A$ implies that $A$ is completely
determined by it values  $a_g=(g\ A\ 1)$, $g\in G$. It is
convenient to think of the form as given by $(x\ A\ y)=\langle x,
Ay\rangle $, where $\langle\ ,\ \rangle $ is the inner product and
$A=\Sigma_{g\in G} a_g g$ is a bounded $G$--equivariant operator on
$\l^2(G)$. Non-degeneracy of $A$ means that $Ax=0$ only if $x=0$.
$A$ is the limit of the group ring elements, and $Ax$ is the limit
of the corresponding linear combinations of $G$--translates of $x$,
i.e.\ $Ax=\lim \Sigma_{g\in G_n} a_g (gx)$, where $G_n$ is some
exhaustion of $G$ by finite sets. It follows that if $x$ belongs
to $G$--invariant isotropic subspace, then $Ax$ belongs to the
closure of this subspace. Thus, we have $\langle Ax,Ax\rangle
=(Ax\ A\ x)=0$ by isotropy and continuity, therefore $x=0$.

The case of general submodule $M \subset \l^2(G)$ reduces to the
above, since the bilinear form $A$ can be extended to $\l^2(G)$,
for example, by taking the orthogonal sum $A \oplus \langle\ ,\
\rangle$ of $A$ on $M$ and the inner product on the orthogonal
complement of $M$.
\end{proof}

\section{Vanishing theorems}

Our main technical results are the following two theorems.

\begin{Theorem}\label{II}
$\III''(2k+1)$ is true for all $k>0$ and $\bq = \bone$.
\end{Theorem}

\begin{proof}
The proof is by induction on $k$. Suppose the theorem is true for
all $m<k$.

Let $v$ be a vertex of $\bdd{2k+1}$. We need to show that the
restriction map $i^*\co \H^k(\bdd{2k+1}) \to \H^k(\bddv{2k+1})$ is
the 0--map.

First let us suppose that $v$ is a vertex of dimension 0, i.e.\ a
vertex of $\Delta^{2k+1}$.

Consider the action of the symmetric group $\Sym{2k+1}$ on
$\Delta^{2k+1}$ which fixes the vertex $v$ and  permutes other
vertices. This action gives a simplicial action of $\Sym{2k+1}$ on
$\bdd{2k+1}$ and therefore, after choosing a base point, lifts to
a cubical action of $\Sym{2k+1}$ on $\Sigma_\bdd{2k+1}$ stabilizing
$\Sigma_{\bddv{2k+1}}$. Let $G'$ be the group of cubical
automorphisms of $\Sigma_\bdd{2k+1}$ generated by this action and
the standard action of $W_\bdd{2k+1}$, and let $G$ be the
orientation-preserving subgroup of $G'$. Similarly, let $G'_v$ be
the group of cubical automorphisms of $\Sigma_{\bddv{2k+1}}$
generated by this action and the standard action of
$W_{\bddv{2k+1}}$, and let $G_v$ be the orientation-preserving
subgroup of $G'_v$.

We claim that, as a Hilbert $G_v$--module
$\h^k(\Sigma_{\bddv{2k+1}})$, is a submodule of $\l^2(G_v)$. Note
that $\bddv{2k+1}$ is naturally isomorphic to $\bdd{2k}$.

Using the inductive assumption and  Lemma~\ref{An1},  we can
remove from $\bdd{2k}$ all odd-dimensional vertices  without
changing the $\ltwo$-cohomology: $\H^*(E_{2k})=\H^*(\bdd{2k})$.
Since the action of $\Sym{2k+1}$ on $\bddv{2k+1}=\bdd{2k}$
preserves the dimension of the vertices, we have isomorphism
$\h^*(G_v\Sigma_{E_{2k}} )=\h^*(\Sigma_{\bdd{2k}})$ as Hilbert
$G_v$--modules.

The complex $E_{2k}$ is spanned by the even-dimensional vertices
of $\bdd{2k}$, which correspond to the proper subsets of vertices
of $\Delta^{2k}$ of odd cardinality. Thus,   the dimension of
$E_{2k}$ is $k-1$, and its top-dimensional simplices are chains
$v_0<v_0v_1v_2<...<v_0...v_{2k-2}$ of length $k$ of distinct
vertices of $\Delta^{2k}$. Therefore $\Sym{2k+1}$ acts
transitively on $(k-1)$--dimensional simplices of $E_{2k}$ and it
follows that $G_v$ acts transitively on $k$--dimensional cubes of
$G_v\Sigma_{E_{2k}}$. Therefore the space of $k$--cochains is a
Hilbert $G_v$--submodule of $\l^2(G_v)$, and the claim follows from
the Hodge decomposition.

We have, by construction, $G_v=\Stab_G(\Sigma_{\bddv{2k+1}})$.
Then the restriction map $i^*\co \h^k(\Sigma_\bdd{2k+1}) \to
\h^k(\Sigma_{\bddv{2k+1}})$ is $G_v$--equivariant and therefore its
image is a $G_v$--invariant subspace of
$\h^k(\Sigma_{\bddv{2k+1}})$. Since $G_v$ acts preserving
orientation, the intersection form is $G_v$--invariant. By
Lemma~\ref{ISO} the image is isotropic, thus by Lemma~\ref{NIS} it
is $0$. Thus, the map $i^*\co \H^k(\bdd{2k+1}) \to
\H^k(\bddv{2k+1})=\h^k(W_\bdd{2k+1}\Sigma_{\bddv{2k+1}})$ is the
$0$--map.

For vertices of the other even dimensions the argument is similar.
If $\dim(v)=2d$, then its link is $\bdd{2d}*\bdd{2k-2d}$. Again,
using Lemma \ref{An1}, we remove, without changing the
$\ltwo$--cohomology, all odd-dimensional vertices from each factor
to obtain $E_{2d}*E_{2k-2d}$. The group
$\Sym{2d+1}\times\Sym{2k-2d+1}$ acts naturally on $\bdd{2k+1}$
fixing the vertex $v$ and stabilizing both the link and
$E_{2d}*E_{2k-2d}$. This action is again transitive on the
top-dimensional simplices of $E_{2d}*E_{2k-2d}$, and the rest of
the argument goes through.

Finally, if $v$ is an odd-dimensional vertex, $\dim(v)=2d+1$, then
we have $\bddv{2k+1}=\bdd{2d+1}*\bdd{2k-2d-1}$. The hypothesis on
$\III''$ and Theorem \ref{t31} and Lemma \ref{l:2.4} imply that
both $\I(2d+1)$ and $\I(2k-2d-1)$ hold. Therefore, by the
K\"unneth formula $\H^k(\bddv{2k+1})=0$ in this case.
\end{proof}

\begin{Theorem}\label{M}
The \GSC\ holds true for $\bdd{2k+1}$:
$$\bt^i(\bdd{2k+1})=0 \text{ for $i> k$ and $\bq \le \bone$.}$$
\end{Theorem}

\begin{proof}
We proceed by induction on $k$.  Using the inductive assumption
and Lemma~\ref{A2k+1}, we can remove all odd-dimensional vertices
without changing the weighted $\ltwoq$--homology above $k$. Thus,
since the remaining part $E_{2k+1}$ is $k$--dimensional, the
problem reduces to showing that $\Ht_{k+1}(E_{2k+1})=0$ for $\bq
\le \bone$. Since $E_{2k+1}$ is $k$--dimensional, the natural map
$\H_{k+1}(E_{2k+1}) \to \Ht_{k+1}(E_{2k+1})$ is injective and the
result follows from the Theorem~\ref{II}.
\end{proof}

Next, we list some consequences. Lemma \ref{An}  implies:

\begin{Theorem}\label{III}
Let $bL$ be the barycentric subdivision of an $(n-1)$--dimensional
simplicial complex $L$. Then $$\bt^i(bL)=0 \text{ for $i >
(n+1)/2$ and $\bq \le \bone$.}$$
\end{Theorem}

Taking $L$ to be a $\GHS^{2k-1}$, we obtain:

\begin{Theorem}\label{IV}
The \GSC\ holds true for the barycentric subdivision of a
$\GHS^{n-1}$ for all even $n$.
\end{Theorem}

 For odd $n$ we obtain a weaker statement:

\begin{Theorem}
Let $bL$ be the barycentric subdivision of a $\GHS^{2k}$. Then
$$\bt^i(bL)=0 \text{ for $i > k+1$ and $\bq \le \bone$.}$$
In particular,
$$\b^i(bL)=0 \text{ for $i \neq k,\ k+1$.}$$
\end{Theorem}

Specializing further, and combining with Theorem \ref{M}, we
obtain:

\begin{Theorem}\label{V}
The \GSC\ holds true for $\bdd{n}$:
$$\bt^i(\bdd{n})=0 \text{ for $i> n/2$ and $\bq \le \bone$,}$$
and, therefore, for the class $\mathcal{JD}$.
\end{Theorem}

Finally, let us mention an application of the above result to a
more analytic object. Let $T_n$ denote the space of all symmetric
tridiagonal $(n+1)\times (n+1)$--matrices with fixed generic
eigenvalues, the so-called Tomei manifold. It is proved in
\cite{To} that $T_n$ is an $n$--dimensional closed aspherical
manifold.

\begin{Theorem}\label{VI}
The Singer Conjecture holds true for Tomei manifolds $T_n$.
\end{Theorem}

\begin{proof}
The space $T_n$ can be identified with a natural finite index
orbifoldal cover of $\Sigma_\bdd{n}/W_\bdd{n}$  \cite{D4}. Thus
$\Sigma_{b\partial \Delta^n}$ is the universal cover of $T^n$, and
the claim follows from the previous theorem.
\end{proof}

\end{document}